\documentclass{amsart} 
\usepackage[mathscr]{eucal}
\usepackage{amsmath,amsfonts}

\newtheorem{theorem}{Theorem}[section]

\newtheorem{lemma}[theorem]{Lemma}

\makeatletter
\@addtoreset{equation}{section}
\makeatother

\newcommand{\CC}{{\mathbb C}}
\newcommand{\NN}{{\mathbb N}}

\newcommand{\ZZ}{{\mathbb Z}}

\newcommand{\cA}{{\mathcal A}}

\newcommand{\cD}{{\mathcal D}}
\newcommand{\cE}{{\mathcal E}}

\newcommand{\cH}{{\mathcal H}}
\newcommand{\cK}{{\mathcal K}}
\newcommand{\cL}{{\mathcal L}}

\newcommand{\cS}{{\mathcal S}}

\newcommand{\cJ}{{\mathcal J}}

\begin{document}
\pagestyle{plain}

\bigskip

\title{Positive definite kernels and lattice paths} 
 \author{T. Constantinescu} \author{Nermine El-Sissi} 

\address{Department of Mathematics \\
  University of Texas at Dallas \\
  Richardson, TX 75083} 
\email{\tt tiberiu@utdallas.edu}
\address{Department of Mathematics \\
  University of Texas at Dallas \\
  Richardson, TX 75083} 
\email{\tt  nae021000@utdallas.edu}

\noindent
\begin{abstract}
We discuss the structure of positive definite kernels 
in terms of operator models. In particular, we introduce 
two models, one of Hessenberg type and another one that we 
call near tridiagonal. These models produce parametrizations 
of the kernels and we describe the combinatorial nature of 
these parametrizations in terms of lattice paths of Dyck and 
Lukasiewicz type.
\end{abstract}
\maketitle

\section{Introduction}
In this paper we are interested in parametrizations and combinatorial
descriptions of positive definite kernels on the set $\NN _0$ 
of non-negative integers. Positive definite kernels are complex valued maps
$K:\NN _0\times \NN _0\rightarrow \CC $ with the 
 property that 
for each $n>0$ and each choice of 
elements $p_1$, $\ldots $, $p_n$ in $\NN _0$ and
complex numbers $\lambda _1$, $\ldots $, $\lambda _n$, we have
\begin{equation}\label{posit}
\sum _{k,j=1}^nK(p_k,p_j)\lambda _j\overline {\lambda }_k\geq 0.
\end{equation}
 A fundamental result of Kolmogorov, \cite{Ko} provides a 
Hilbert space interpretation of positive definite kernels as Gram 
kernels, that is, there exists a Hilbert space $\cH $ and elements
$v(n)\in \cH $, $n\geq 0$, such that 
\begin{equation}\label{kolmo}
K(i,j)=\langle v(j), v(i)\rangle .
\end{equation}

Two of the best known examples of positive definite kernels
are those of Toeplitz type, for which
$K(i+l,j+l)=K(i,j)$, $i,j,l\in \NN _0,$
and those of Hankel type, for which 
$K(i,j+l)=K(i+l,j)$, $i,j,l\in \NN _0.$
In both these cases the representation \eqref{kolmo}
can be improved by some more specialized descriptions that might
be called operator models of the kernels. Thus, if $K$ is Toeplitz
(for simplicity we assume $K(0,0)=1$ and all the inequalities 
in \eqref{posit} are strict) then there exists an isometric operator $W$,
written in upper Hessenberg form, such that 
\begin{equation}\label{toeplitz}
K(i,j)=e_0W^{j-i}e^*_0,
\end{equation} 
where $e_0=\left[
\begin{array}{ccc}
1 & 0 & \ldots 
\end{array}
\right]$. Likewise, if $K$ is Hankel (with same simplifications as above)
then there exists a symmetric operator $J$, written in tridiagonal
form, such that 
\begin{equation}\label{hankel}
K(i,j)=e_0J^{i+j}e^*_0.
\end{equation} 

Our goal is to extend both models \eqref{toeplitz}
and \eqref{hankel} to arbitrary positive definite kernels on $\NN _0$
without Toeplitz or Hankel assumptions. 
These models will produce parametrizations of the kernels and 
we will give combinatorial descriptions of these parametrizations in terms of 
lattice paths.

\section{Isometric Hessenberg models and Dyck paths}
In this section we show that any positive definite 
kernel on $\NN _0$ has a Hessenberg model and then we show
how to relate this model to the 
set of Dyck paths. In order to simplify the notation we consider 
only positive definite kernels for which all the inequalities in 
\eqref{posit} are strict, when we say that the kernel 
is stricly positive definite, and we also assume 
$K(l,l)=1$ for all $l\geq 0$. Both these assumptions can be easily
removed. In addition, all our considerations can be easily adapted to kernels
$K:\NN _0\times \NN _0\rightarrow \cL (\cE )$, where $\cL (\cE )$ denotes the 
set of linear bounded operators on the Hilbert space $\cE $.

We now introduce the elements necesary in the presentation of the results. 
For a complex number $\gamma $
with $|\gamma |\leq 1$ we define its {\it defect} by 
$d_{\gamma }=(1-|\gamma |^2)^{1/2}$ and its  
{\it Julia matrix } by 
$$J(\gamma )=\left[\begin{array}{cc}
\gamma & d_{\gamma } \\
d_{\gamma } & -\overline{\gamma }
\end{array}
\right].
$$
We note that the Julia matrix is unitary and this construction can be 
extended to
certain families of complex numbers as follows. Let 
$\Gamma =\{\gamma _{k,j}\}_{0\leq k<j}$ be a family of 
complex numbers such that 
$|\gamma _{k,j}|<1$ for all
$k<j$. For simplicity we 
will write $d_{k,j}$ instead of $d_{\gamma _{k,j}}$.
We can now describe the Hessenberg model. First we define for $k<j$,
$$V_{k,j}(\Gamma )=
\left(J(\gamma _{k,k+1})\oplus I_{j-k-1}\right)
\left(I_1\oplus J(\gamma _{k,k+2})\oplus I_{j-k-2}\right)\ldots
\left( I_{j-k-1}\oplus J(\gamma _{k,j})\right),$$
where $I_l$ denotes the $l\times l$ identity matrix. Then we introduce 
the operators $W_k(\Gamma )$ on the Hilbert space $l^2(\NN _0)$
of square-summable sequences by the formula:
$$W_k(\Gamma )=s-\lim _{j\rightarrow \infty }(V_{k,j}(\Gamma )\oplus 0),
\quad k\geq 0,$$
where $s-\lim$ denotes the strong operator limit. It is easily seen that 
each $W_k(\Gamma )$ is an isometry with upper Hessenberg
matrix with respect to the standard basis of $l^2(\NN _0)$, that
is, if $(W_k(\Gamma ))_{i,j}$ denotes the $(i,j)$ entry of $W_k(\Gamma )$, 
then $(W_k(\Gamma ))_{i,j}=0$ for $i\geq j+1$. It is also useful 
to consider the unitary matrices $U_{k,j}(\Gamma )$ defined recursively 
by $U_{k,k}(\Gamma )=I_1$ and for $k<j$,
$$U_{k,j}(\Gamma )=V_{k,j}(\Gamma )\left(U_{k+1,j}(\Gamma )\oplus I_1\right).
$$
We can prove now the existence of an isometric 
Hessenberg model for any strictly positive definite
kernel on $\NN _0$.
\begin{theorem}\label{Hesse}
If $K$ is a strictly positive definite kernel on $\NN _0$ with 
$K(l,l)=1$ for all $l\geq 0$, then there exists a family $\{W_k\}_{k\geq 0}$
of isometric Hessenberg operators such that for $j>i$, 
\begin{equation}\label{dilat}
K(i,j)=e_0W_i(\Gamma )W_{i+1}(\Gamma )\ldots W_{j-1}(\Gamma )e^*_0.
\end{equation}
\end{theorem}
\begin{proof}
This is just a reformulation of Theorem ~2.3 in \cite{Co1}
(see also \cite{Co}, Chapter 1). Thus, by Theorem ~1.3 in \cite{Co1}, 
there exists a uniquely determined family 
$\Gamma =\{\gamma _{k,j}\}_{0\leq k<j}$ of complex numbers 
such that $K(i,j)=\left(U_{i,j}(\Gamma )\right)_{0,0}$
for $i<j$. Then it is easily seen from the definitions that 
$$ 
e_0W_i(\Gamma )W_{i+1}(\Gamma )\ldots W_{j-1}(\Gamma )e^*_0=
\left(U_{i,j}(\Gamma )\right)_{0,0}.
$$
\end{proof}

When $K$ is a Toeplitz kernel, then \eqref{dilat}
reduces to \eqref{toeplitz} and the parameters $\gamma _{k,j}$
satisfy $\gamma _{i+l,j+l}=\gamma _{i,j}$ for 
$i<j$, $l\geq 1$. The numbers $\gamma _n=\gamma _{0,n}$, $n\geq 1$, 
are called the {\it Szeg\"o parameters} of $K$ (other names, like
Schur parameters, reflection coefficients, or Verblunski parameters
are currently used in the literature), and they play a central role
in the theory of orthogonal polynomials on the unit circle and its many 
applications, \cite{Sz} and, for a recent account, \cite{Si}
(which also contains a detailed discussion of the Hessenberg model 
in the Toeplitz case). 

\smallskip
Next we explain the connection between Hessenberg models and Dyck paths.
A Dyck path of length $2k$ is a path in the 
positive quadrant of the lattice $\ZZ ^2$ which starts at $(0,0)$, 
ends at $(2k,0)$, and consists of rise steps $\nearrow $
and fall steps $\searrow $ (see Figure ~1). For more 
information on Dyck paths and their combinatorics, see \cite{St}.  

\begin{figure}[h]
\setlength{\unitlength}{3000sp}%
\begingroup\makeatletter\ifx\SetFigFont\undefined%
\gdef\SetFigFont#1#2#3#4#5{%
  \reset@font\fontsize{#1}{#2pt}%
  \fontfamily{#3}\fontseries{#4}\fontshape{#5}%
  \selectfont}%
\fi\endgroup%
\begin{picture}(2424,1224)(289,-673)
{ \thinlines
\multiput(301,239)(8.98876,0.00000){268}{\makebox(1.6667,11.6667){\SetFigFont{5}{6}{\rmdefault}{\mddefault}{\updefault}.}}

\multiput(301,-61)(8.98876,0.00000){268}{\makebox(1.6667,11.6667){\SetFigFont{5}{6}{\rmdefault}{\mddefault}{\updefault}.}}

\multiput(301,-361)(8.98876,0.00000){268}{\makebox(1.6667,11.6667){\SetFigFont{5}{6}{\rmdefault}{\mddefault}{\updefault}.}}

\multiput(601,-661)(0.00000,9.02256){134}{\makebox(1.6667,11.6667){\SetFigFont{5}{6}{\rmdefault}{\mddefault}{\updefault}.}}

\multiput(901,-661)(0.00000,9.02256){134}{\makebox(1.6667,11.6667){\SetFigFont{5}{6}{\rmdefault}{\mddefault}{\updefault}.}}

\multiput(1201,-661)(0.00000,9.02256){134}{\makebox(1.6667,11.6667){\SetFigFont{5}{6}{\rmdefault}{\mddefault}{\updefault}.}}

\multiput(1501,-661)(0.00000,9.02256){134}{\makebox(1.6667,11.6667){\SetFigFont{5}{6}{\rmdefault}{\mddefault}{\updefault}.}}

\multiput(1801,-661)(0.00000,9.02256){134}{\makebox(1.6667,11.6667){\SetFigFont{5}{6}{\rmdefault}{\mddefault}{\updefault}.}}

\multiput(2101,-661)(0.00000,9.02256){134}{\makebox(1.6667,11.6667){\SetFigFont{5}{6}{\rmdefault}{\mddefault}{\updefault}.}}

\multiput(2401,-661)(0.00000,9.02256){134}{\makebox(1.6667,11.6667){\SetFigFont{5}{6}{\rmdefault}{\mddefault}{\updefault}.}}

\put(301,-661){\framebox(2400,1200){}}

\put(301,-661){\line( 1, 1){600}}
\put(901,-61){\line( 1,-1){300}}
\put(1201,-361){\line( 1, 1){600}}
\put(1801,239){\line( 1,-1){900}}
}%

\end{picture}

\caption{\mbox{ A Dyck path of length $8$ and height $3$}}
\end{figure}

Let $\cD _k$ 
be the set of Dyck paths of length $2k$ and 
let $\cA _k$ be the set of points $(l,q)$, $q>0$, 
with the property that there exists ${\bf p}\in \cD _k$ with 
$(l,q)\in {\bf p}$. It is seen that 
$$\cA _k=\left\{(j+i,j-i)\mid 0\leq i<j\leq k\right\}.$$
Also, we notice that if  ${\bf p}\in \cD _k$ and 
$x=(l,q)\in {\bf p}$, then there are only four types of behaviour of $
{\bf p}$
about $x$: (I) a rise step followed by a fall step; 
(II) a fall step followed by a rise step; (III)
two consecutive rise steps; (IV) two consecutive fall steps (see
Figure ~2). 

\begin{figure}[h]
\setlength{\unitlength}{3000sp}%
\begingroup\makeatletter\ifx\SetFigFont\undefined%
\gdef\SetFigFont#1#2#3#4#5{%
  \reset@font\fontsize{#1}{#2pt}%
  \fontfamily{#3}\fontseries{#4}\fontshape{#5}%
  \selectfont}%
\fi\endgroup%
\begin{picture}(4224,966)(589,-1015)
{ \thinlines
{\put(601,-661){\line( 1, 1){300}}
}%
{\put(901,-361){\line( 1,-1){300}}
}%
{\put(1801,-361){\line( 1,-1){300}}
}%
{\put(2101,-661){\line( 1, 1){300}}
}%
{\put(3001,-661){\line( 1, 1){600}}
}%
{\put(4201,-61){\line( 1,-1){600}}
}%
\put(826,-961){\makebox(0,0)[lb]{\smash{{\SetFigFont{12}{14.4}{\rmdefault}{\mddefault}{\updefault}{$(I)$}%
}}}}
\put(2026,-961){\makebox(0,0)[lb]{\smash{{\SetFigFont{12}{14.4}{\rmdefault}{\mddefault}{\updefault}{$(II)$}%
}}}}
\put(3226,-961){\makebox(0,0)[lb]{\smash{{\SetFigFont{12}{14.4}{\rmdefault}{\mddefault}{\updefault}{$(III)$}%
}}}}
\put(4276,-961){\makebox(0,0)[lb]{\smash{{\SetFigFont{12}{14.4}{\rmdefault}{\mddefault}{\updefault}{$(IV)$}%
}}}}
\put(850, -600){$x$}
\put(2050, -500){$x$}
\put(3320, -500){$x$}
\put(4300, -500){$x$}

}
\end{picture}

\caption{\mbox{ Behaviour of a Dyck path about a vertex $x\in \cA _k$}}
\end{figure}

Consequently, 
for each pair $i,j$ with $0\leq i<j\leq k$ we define the function
$a_{i,j}:\cD _k\rightarrow \CC $, 
$$a_{i,j}({\bf p})=\left\{\begin{array}{ccl}
1 &\mbox{if} & x=(j+i,j-i)\notin {\bf p}; \\
\gamma _{i,j} &\mbox{if} & x=(j+i,j-i)\in {\bf p} \quad \mbox{and (I) holds};\\
-\overline{\gamma }_{i,j} &\mbox{if} & x=(j+i,j-i)\in {\bf p} \quad \mbox{and (II) holds};\\
d_{i,j} &\mbox{if} & x=(j+i,j-i)\in {\bf p} \quad \mbox{and
either (III) or (IV) holds}.
\end{array}
\right.
$$
Let ${\bf p}$ be a Dyck path in $\cD _k$ with the property
that $(2l,0)\in {\bf p}$.
The restriction of ${\bf p}$ from $(2l,0)$ to $(2k,0)$ is called
a Dyck subpath starting at $(2l,0)$ in $\cD _k$. The set of all 
possible Dyck subpaths starting at $(2l,0)$ in $\cD _k$ is denoted by  
$\cD ^l_k$ and there exists a bijection 
between $\cD ^l_k$ and $\cD _{k-l}$. This implies that the number of elements
in $\cD ^l_k$ is given by the Catalan number 
$C_{k-l}=\displaystyle\frac{1}{k-l+1}\left(
\begin{array}{c}
2(k-l) \\
k-l
\end{array}
\right);
$
also, $\cD ^0_k=\cD _k$.
If ${\bf q}\in \cD ^l_k$ then there could be many Dyck paths 
whose restrictions at $(2l,0)$ coincide with ${\bf q}$, however 
if ${\bf p}_1$ and ${\bf p}_2$ are two such Dyck paths then 
$a_{i,j}({\bf p}_1)=a_{i,j}({\bf p}_2)$ for $j+i>2l$. We will write 
$a_{i,j}({\bf q})$ in order to denote this common value. 

We now describe the structure of the strictly positive definite kernels
on $\NN _0$.
\begin{theorem}\label{suma}
The kernel $K$ on $\NN _0$ with $K(l,l)=1$ for all $l$
is strictly positive definite
if and only if there is a family
$\{\gamma _{k,j}\}_{0\leq k\leq j}$ of 
complex numbers, $|\gamma _{k,j}|<1$ for all
$k<j$, such that 
\begin{equation}\label{fsuma}
K(l,m)=\sum _{{\bf q}\in \cD ^l_m}\prod _{l\leq i<j\leq m}a_{i,j}({\bf q}).
\end{equation}
\end{theorem} 
\begin{proof}
Half of this result was proved in \cite{BCE}, but we give some details
here for completeness. Assume that $K$ is strictly positive definite.
By Theorem ~1.3 and Theorem ~2.3 in \cite{Co1} there exists a uniquely 
determined family 
$\Gamma =\{\gamma _{k,j}\}_{0\leq k<j}$ of complex numbers such that 
$K(l,m)=(U_{l,m}(\Gamma ))_{0,0}$, the $(0,0)$ entry of the 
matrix $U_{l,m}(\Gamma )$. It is convenient to visualize this 
relation by means 
of a so-called transmission line, as showed in Figure ~3 for 
$K(0,2)$ and $K(0,3)$.
\begin{figure}[h]
\setlength{\unitlength}{3000sp}%
\begingroup\makeatletter\ifx\SetFigFont\undefined%
\gdef\SetFigFont#1#2#3#4#5{%
  \reset@font\fontsize{#1}{#2pt}%
  \fontfamily{#3}\fontseries{#4}\fontshape{#5}%
  \selectfont}%
\fi\endgroup%
\begin{picture}(5028,2364)(187,-1615)
{ \thinlines

\put(301,164){\circle{212}}
\put(2101,164){\circle{212}}
\put(5101,164){\circle{212}}

\put(601,-361){\framebox(600,600){}}
\put(2401,-361){\framebox(600,600){}}
\put(4201,-361){\framebox(600,600){}}
\put(1501,-811){\framebox(600,600){}}
\put(3301,-811){\framebox(600,600){}}
\put(2401,-1261){\framebox(600,600){}}

\put(301,164){\line( 1, 0){4800}}
\put(301,-286){\line( 1, 0){4800}}
\put(301,-736){\line( 1, 0){4800}}
\put(301,-1186){\line( 1, 0){4800}}
\put(676,-286){\line( 1, 1){450}}
\put(676,164){\line( 1,-1){450}}
\put(2476,164){\line( 1,-1){450}}
\put(4276,164){\line( 1,-1){450}}
\put(2476,-286){\line( 1, 1){450}}
\put(4276,-286){\line( 1, 1){450}}
\put(1576,-736){\line( 1, 1){450}}
\put(3376,-736){\line( 1, 1){450}}
\put(2476,-1186){\line( 1, 1){450}}
\put(1576,-286){\line( 1,-1){450}}
\put(3376,-286){\line( 1,-1){450}}
\put(2476,-736){\line( 1,-1){450}}
\put(301, 89){\line( 0, 1){450}}
\put(2101, 89){\line( 0, 1){450}}
\put(5101, 89){\line( 0, 1){450}}
\put(226,164){\line( 1, 0){ 75}}
\put(5101,164){\line( 1, 0){ 75}}
\put(301,539){\vector( 0,-1){300}}
\put(2101,539){\vector( 0,-1){300}}
\put(5101,239){\vector( 0, 1){300}}
}%
\put(676,389){$J(\gamma _{2,3})$}%

\put(2476,389){$J(\gamma _{1,2})$}%

\put(4276,389){$J(\gamma _{0,1})$}%

\put(1576,-61){$J(\gamma _{1,3})$}%

\put(3376,-61){$J(\gamma _{0,2})$}%

\put(2476,-1561){$J(\gamma _{0,3})$}%

\put(226,614){$A$}%

\put(2026,614){$C$}%

\put(5026,614){$B$}%

\end{picture}

\caption{\mbox{ Transmission line for $K(0,3)$ }}
\end{figure}

\noindent
Thus, if $1$ is the input at $A$ then at $B$
we read off the expression of $K(0,3)$ in terms of the 
parameters $\gamma _{0,1}$, $\gamma _{0,2}$,
$\gamma _{0,3}$, $\gamma _{1,2}$, $\gamma _{1,3}$, $\gamma _{2,3}$
and their defects,
$$\begin{array}{rl}
K(0,3)=&\!\!\!\gamma _{0,1}\gamma _{1,2}\gamma _{2,3}
+\gamma _{0,1}d_{1,2}\gamma _{1,3}d_{2,3}+
d_{0,1}\gamma _{0,2}d_{1,2}\gamma _{2,3} \\
  & -d_{0,1}\gamma _{0,2}\overline{\gamma }_{1,2}\gamma _{1,3}d_{2,3}+
d_{0,1}d_{0,2}\gamma _{0,3}d_{1,3}d_{2,3}.
\end{array}
$$
Likewise, if the input at $C$ is $1$, 
then the output at $B$ is now the expression of $K(0,2)$,
$$K(0,2)=\gamma _{0,1}\gamma _{1,2}+d_{0,1}\gamma _{0,2}d_{1,2}
$$
(for more details see \cite{Co}).
Each path in the transmission line contributes an additive term in 
$K(l,m)$. Going from a path in the transmission line to a Dyck path is easy, 
each box associated with a Julia matrix corresponds to a point 
in $\cA _k$, see Figure ~4. 
It is also clear that each additive
term in 
$K(l,m)$ is given by $\prod _{l\leq i<j\leq m}a_{i,j}({\bf q})$
for some ${\bf q}\in \cD ^l_m$. This gives \eqref{fsuma}.

Conversely, given a family $\Gamma =\{\gamma _{k,j}\}_{0\leq k<j}$ of complex 
numbers with $|\gamma _{k,j}|<1$ for all $k<j$, we define
$$K(l,m)=\sum _{{\bf q}\in \cD ^l_m}\prod _{l\leq i<j\leq m}a_{i,j}({\bf q}).
$$
By the first part of the proof, this gives $K(l,m)=(U_{l,m}(\Gamma ))_{0,0}$, 
and it remains to show that $K$ is a strictly 
positive definite kernel on $\NN _0$. 
By Theorem ~2.1, 
$$K(i,j)=e_0W_i(\Gamma )W_{i+1}(\Gamma )\ldots W_{j-1}(\Gamma )e^*_0
$$
for $i<j$. This relation 
implies
that $K$ is a positive definite kernel. 
Also, by Proposition ~1.7 in  \cite{Co1}, 
$$\det \left[K(l,m)\right]_{l,m=0}^n=\prod _{0\leq i<j\leq n}d^2_{i,j}>0,$$
so that $K$ is a strictly positive definite kernel on $\NN _0$.
\end{proof}

\begin{figure}[h]
\setlength{\unitlength}{3000sp}%
\begingroup\makeatletter\ifx\SetFigFont\undefined%
\gdef\SetFigFont#1#2#3#4#5{%
  \reset@font\fontsize{#1}{#2pt}%
  \fontfamily{#3}\fontseries{#4}\fontshape{#5}%
  \selectfont}%
\fi\endgroup%
\begin{picture}(6405,1525)(218,-828)
\thinlines
{\put(2814,-511){\oval(976,602)[bl]}
\put(2814,-511){\oval(974,602)[br]}
\put(3301,-511){\vector( 0, 1){0}}
}%
{\put(4904,-511){\oval(1106,602)[bl]}
\put(4904,-586){\oval(994,452)[br]}
\put(5401,-586){\vector( 0, 1){0}}
}%
{\thinlines
\put(301,614){\circle{150}}
}%
{\put(1351,614){\circle{150}}
}%
{\put(2401,614){\circle{150}}
}%
{\put(3001,614){\circle{150}}
}%
{\put(3751,614){\circle{150}}
}%
{\put(4501,614){\circle{150}}
}%
{\put(3751,-361){\circle{212}}
}%
{\put(3376, 89){\circle{212}}
}%
{\put(4095, 89){\circle{212}}
}%
{\put(751, 89){\circle{618}}
}%
{\put(1351,-202){\circle{618}}
}%
{\put(1949, 89){\circle{618}}
}%
\thicklines
{\put(301,539){\line( 0,-1){300}}
\put(301,239){\line( 1, 0){300}}
\put(601,239){\line( 1,-1){300}}
\put(901,-61){\line( 1, 0){900}}
\put(1801,-61){\line( 1, 1){300}}
\put(2101,239){\line( 1, 0){300}}
\put(2401,239){\line( 0, 1){300}}
}%
\thinlines
{\put(601,239){\line( 1, 0){1500}}
}%
{\put(301,-61){\line( 1, 0){2100}}
}%
{\put(301,-361){\line( 1, 0){2100}}
}%
{\put(601,-61){\line( 1, 1){300}}
}%
{\put(1801,239){\line( 1,-1){300}}
}%
{\put(1201,-361){\line( 1, 1){300}}
}%
{\put(1201,-61){\line( 1,-1){300}}
}%
{\put(526,-136){\framebox(450,450){}}
}%
{\put(1126,-436){\framebox(450,450){}}
}%
{\put(1726,-136){\framebox(450,450){}}
}%
{\put(1351,539){\line( 0,-1){300}}
}%
\thicklines
{\put(3001,614){\line( 3,-4){738}}
\put(3751,-361){\line( 3, 4){738}}
}%
\thinlines
{\put(5401,-361){\line( 1, 0){1200}}
\put(6601,-361){\line( 0, 1){600}}
\put(6601,239){\line(-1, 0){1200}}
\put(5401,239){\line( 0,-1){600}}
}%
{\put(5401,-61){\line( 1, 0){1200}}
}%
{\put(5701,239){\line( 0,-1){600}}
}%
{\put(6001,239){\line( 0,-1){600}}
}%
{\put(6301,239){\line( 0,-1){600}}
}%
\thicklines
{\put(5401,-361){\line( 1, 1){600}}
\put(6001,239){\line( 1,-1){600}}
}%

\end{picture}%

\caption{\mbox{ From a path in a transmission line to a Dyck path }}
\end{figure}

\bigskip
\noindent
{\bf Remarks}
$(a)$ It is quite simple to remove the two restrictions on $K$ 
considered in Theorem ~\ref{suma}. First, formula \eqref{fsuma} still
provides a one-to-one correspondence between the set of positive definite 
kernels on $\NN _0$ with $K(l,l)=1$ for all $l$ and the set
$\cS $ of families $\{\gamma _{k,j}\}_{0\leq k<j}$ of 
complex numbers with the properties: $|\gamma _{k,j}|\leq 1$ for 
all $k<j$; if $|\gamma _{k,j}|=1$ for some pair $(k,j)$, then 
$\gamma _{l,j}=0$ for $l<k$ and $\gamma _{k,m}=0$ for $m>j$.

\smallskip
\noindent
$(b)$ If we remove the assumption that $K(l,l)=1$ for all $l$, 
then the diagonal elements 
$K(l,l)$ of the kernel $K$ could be considered as parameters and there
is a one-to-one 
correspondence between the positive definite kernels on $\NN _0$ and the
set $\cS _+$ of pairs $\left(\{f_l\}_{l\geq 0},
\{\gamma _{k,j}\}_{0\leq k<j}\right)$, where $f_l\geq 0$ for all 
$l\geq 0$ and $\{\gamma _{k,j}\}_{0\leq k<j}$ is an element of $\cS $ 
with the additional property that if $f_l=0$ for some $l\geq 0$, then 
$\gamma _{k,l}=0$ and $\gamma _{l,m}=0$ for $k<l$ and
$m>l$. Formula \eqref{fsuma} has to be replaced in this case with:
\begin{equation}\label{ffsuma}
K(l,m)=
f^{1/2}_lf^{1/2}_m
\sum _{{\bf q}\in \cD ^l_m}\prod _{l\leq i<j\leq m}a_{i,j}({\bf q}).
\end{equation}

\smallskip
\noindent
$(c)$ In case 
$K$ is a Toeplitz kernel we noticed already that 
$\gamma _{i+l,j+l}=\gamma _{i,j}$ for $i<j$, $l\geq 1$
and we denoted $\gamma _n=\gamma _{o,n}$, $n\geq 1$. We conclude
that $a_{i+l,j+l}=a_{i,j}$ for $i<j$, $l\geq 1$ and formula
\eqref{fsuma} reduces to
\begin{equation}\label{tsuma}
K(0,n)=\sum _{{\bf p}\in \cD _n}\prod _{0\leq i<j\leq n}a_{i,j}({\bf p}).
\end{equation}
We can compare this result 
with a classical formula of Verblunsky, according to which
there exists a polynomial 
$V^{(n)}=V^{(n)}(\gamma _1,\ldots \gamma _{n-1};d_1,\ldots ,d_{n-1})$
with integer coeffcients so that
$$K(0,n)=
\gamma _n\prod _{k=1}^{n-1}d_k^2+
V^{(n)}$$ 
(see \cite{Si}, in particular, pg. 60-61, 
for a comprehensive discussion of this formula).
We see that the term $\gamma _n\prod _{k=1}^{n-1}d_k^2$
corresponds to the path ${\bf p}_0$ made 
of $n$ consecutive rise steps followed by $n$ consecutive fall
steps. Consequently, we deduce from \eqref{tsuma}
that 
$$
V^{(n)}=
\sum _{{\bf p}\in \cD _n-\{{\bf p}_0\}}
\prod _{0\leq i<j\leq n}a_{i,j}({\bf p}),
$$ 
an explicit formula that explains some of the features of 
$V^{(n)}$.
\quad \quad \quad \quad \quad \quad \quad \quad \quad
$\Box $

\section{Near tridiagonal models}
In this section we show that positive definite kernels do not 
have tridiagonal models. Instead we introduce a near tridiagonal 
model and then we show how this model
is related to the set of Lukasiewicz paths.  
Again, in order to simplify the notation we consider 
only strictly positive definite kernels $K$
and we assume $K(0,0)=1$.
We denote 
by $\cD \subset l^2(\NN _0)$ the vector space generated by the standard basis 
of $l^2(\NN _0)$ and we call {\it tridiagonal model} of $K$ a family
$\{J_n\}_{n\geq 0}$ of tridiagonal operators 
(not necessarely bounded),
$$J_n=\left[
\begin{array}{cccc}
b_0(n) & c_1(n) & 0 & \\
a_1(n) & b_1(n) & c_2(n) & \\
0 & a_2(n) & b_2(n) & \ddots \\
 & 0 & \ddots & \ddots 
\end{array}
\right],
$$
such that $J_0=I$ and 
\begin{equation}\label{haha}
K(i,j)=e_0J_1^*\ldots J^*_iJ_j\ldots J_1e^*_0, \quad i,j\geq 0.
\end{equation}
Also, in analogy with the Hankel case, we ask $a_k(n)>0$, $k,n\geq 1$.
Since each $J_n$ is tridiagonal, $J_n\cD\subset \cD $ so \eqref{haha}
makes sense. However we have the following result.
\begin{theorem}\label{no}
There are strictly positive definite kernels with no tridiagonal model.
\end{theorem}
\begin{proof}
We consider a strictly positive definite kernel $K$ with 
$$K(0,1)=K(0,2)=K(1,2)=0, \quad K(0,3)\ne 0$$
(it is easy to construct such a kernel using, for instance, 
Theorem ~\ref{suma}). Let $\{J_n\}_{n\geq 1}$ be a tridiagonal model
of $K$, then we deduce that 
$$b_0(1)=c_1(2)=b_1(2)=0,$$
which implies 
$$\begin{array}{rcl}
K(0,3) & = & e_0J_3J_2J_1e^*_0\\
       & = & b_0(3)\left(b_0(2)b_0(1)+c_1(2)a_1(1)\right)
+c_1(3)\left(a_1(2)b_0(1)+b_1(2)a_1(1)\right)=0,
\end{array}
$$
a contradiction showing that $K$ has no tridiagonal model.
\end{proof}

We are now trying to find a model as close as possible of being tridiagonal, 
which should reduce to \eqref{hankel} in case the kernel $K$ is 
Hankel. Thus, we consider operators (not necessarely bounded) with 
matrix still of Hessenberg form
\begin{equation}\label{hump}
J_n=\left[
\begin{array}{cccc}
b_0(1) & c_{0,1}(n) & c_{0,2}(n) & \\
a_1(1) & b_1(2) & c_{1,2}(n) & \\
0 & a_2(2) & b_2(3) & \ddots \\
 & 0 & \ddots & \ddots 
\end{array}
\right], \quad n\geq 1,
\end{equation}
with respect to the standard basis of $l^2(\NN _0)$, and with the 
additional conditions:
\begin{equation}\label{conditii}
\begin{array}{l}
a_k(k)>0, \quad k\geq 1, \\
c_{i,j}(1)=0, \quad j\geq 2, \quad 0\leq i<j-1, \\
c_{i,j}(n)=0, \quad j\geq n, \quad 0\leq i<j-1,\\
c_{k-1,k}(n)=a_k(k), \quad k\geq n,            \\
c_{i,j}(n)=c_{i,j}(n-1), \quad j<n-1, \quad 0\leq i<j-1.
\end{array}
\end{equation}
We see that $J_n\cD \subset \cD $ for each $n\geq 1$. We call such a family 
$\{J_n\}_{n\geq 1}$ a {\it near tridiagonal model} 
of the kernel $K$ provided that 
\begin{equation}\label{balabusta}
K(i,j)=e_0J_1^*\ldots J^*_iJ_j\ldots J_1e^*_0, \quad i,j\geq 0.
\end{equation}
\begin{theorem}\label{yes}
Any strictly positive definite kernel has a near tridiagonal model.
\end{theorem}
\begin{proof}
Let $K$ be a kernel and $K_n=\left[K(i,j)\right]_{0\leq i,j\leq n}$.
Then $K$ is strictly positive definite if and only if $K_n>0$, $n\geq 1$
(as already mentioned we can assume without loss of generality that 
$K(0,0)=1$). Let $D_n=\left[d_{i,j}(n)\right]_{0\leq i,j\leq n}$
be the upper triangular Cholesky factor of $K_n$, therefore
$K_n=D^*_nD_n$ and $d_{i,i}(n)>0$. The uniqueness of the Cholesky
factor implies that 
$$D_{n+1}=\left[
\begin{array}{cc}
D_n & l_{n+1} \\
0 & d_{n+1,n+1}(n+1)
\end{array}
\right],
$$
and $d_{n,n}(k+1)=d_{n,n}(k)$ for $k\geq n$, so we can 
drop the label $n$ in $d_{i,j}(n)$.
We now construct the near tridiagonal model of $K$.  Thus we prove 
by induction on $n$ that there exist numbers $a_k(k)$, $b_{k-1}(k)$, 
$c_{i,k-1}(k)$, $0\leq i\leq k-1$, $k\leq n$, such that 
\eqref{conditii} holds and 
 \begin{equation}\label{trunc}
D_n=\left[
\begin{array}{ccccc}
1 & J_{2,1} & J_{3,2}J_{2,1} &  \ldots & J_{n+1,n}\ldots J_{2,1}      \\ 
0_n  & 0_{n-1} & 0_{n-2}     &  \ldots &  0_0
\end{array}
\right],
\end{equation}
where 
$$J_{2,1}=\left[
\begin{array}{c}
b_0(1) \\
a_1(1) 
\end{array}
\right], \quad 
J_{k+1,k}=\left[
\begin{array}{cccc}
b_0(1) & c_{0,1}(2) & \ldots & c_{0,k-1}(k) \\
a_1(1) & b_1(2)     &        &              \\
0      & a_2(2)     & \ddots &              \\
       &  \ddots    & \ddots &              \\
0       &            &   0     & a_{k}(k) 
\end{array}
\right],
$$
and $0_k$ denotes the column with $k$ zero entries.

For $n=1$ we define 
$$b_0(1)=d_{0,1}\quad \mbox{and} \quad a_1(1)=d_{1,1}>0,$$
so that $D_1=\left[
\begin{array}{cc}
1 & \\
  & J_{2,1} \\
0 & 
\end{array}
\right]$. Assume the statement is true up to $n$. We determine 
the numbers $a_{n+1}(n+1)>0$, $b_n(n+1)$, and $c_{k,n}(n+1)$, 
$k=0,\ldots n, n-1$, such that 
$$\left[
\begin{array}{c}
l_{n+1} \\
d_{n+1,n+1}
\end{array}
\right]=J_{n+2,n+1}J_{n+1,n}\ldots J_{2,1}.$$
By the induction hypothesis
$$J_{n+1,n}\ldots J_{2,1}=\left[
\begin{array}{c}
l_{n} \\
d_{n,n}
\end{array}
\right],
$$
so that we must have
\begin{equation}\label{key}
\left[
\begin{array}{c}
l_{n+1} \\
d_{n+1,n+1}
\end{array}
\right]=\left[
\begin{array}{c}
x_0+c_{0,n}(n+1)d_{n,n} \\
x_1+c_{1,n}(n+1)d_{n,n} \\
\vdots \\
x_n+b_n(n+1)d_{n,n}  \\
a_{n+1}(n+1)d_{n,n}
\end{array}
\right],
\end{equation}
where $x_0$, $\ldots $, $x_n$ are numbers uniquely determined 
by $a_{k}(k)>0$, $b_{k-1}(k)$, and $c_{i,l}(k)$, 
$i<l<k\leq n$. Since $d_{n,n}>0$, \eqref{key}
uniquely determine the numbers 
$c_{0,n}(n+1)$, $\ldots $, 
$c_{n-1,n}(n+1)$, $b_n(n+1)$ and $a_{n+1}(n+1)=
\frac{d_{n+1,n+1}}{d_{n,n}}>0$
such that \eqref{trunc} holds for $n+1$. 

Now we can use all the numbers 
$a_k(k)$, $b_{k-1}(k)$, $c_{i,l}(k)$ 
in order to define $J_n$ by \eqref{hump} and then \eqref{trunc}
shows that $\{J_n\}_{n\geq 1}$ is a near tridiagonal model of $K$.
\end{proof}
We notice that the label $n$ of the numbers 
$a_n(n)$, $b_{n-1}(n)$, $c_{i,l}(n)$ is superfluous
due to the conditions \eqref{conditii}. We used it in order to have a uniform 
definition of $J_n$ in \eqref{hump} but we will drop it from now on.
The proof of Theorem ~\ref{yes}
gives a one-to-one correspondence between the set of strictly 
positive definite kernels on $\NN _0$ with $K(0,0)=1$ and the 
set $\cJ $ 
of families of numbers $\{a_k, b_{k-1}, c_{i,l}
\mid k\geq 1,\quad 0\leq i<l\}$.
We will call these numbers the {\it Jacobi parameters} of $K$.
In addition, we can easily characterize the strictly positive definite 
Hankel kernels by the additional 
conditions on the Jacobi parameters:
\begin{equation}\label{carhank}
\begin{array}{l}
c_{i,l}=0,\quad l>i+1,\\
c_{k,k+1}=a_{k+1}, \quad k\geq 0.
\end{array}
\end{equation}
In this case the near tridiagonal model reduces to \eqref{hankel}.

The next task is to establish an explicit formula for the 
Cholesky factors $D_n$ in terms of the Jacobi parameters.
First, we obtain a recursive relation for $D_n$.
\begin{lemma}\label{luna}
For $n\geq 1$, 
$$D_n=F_n(1\oplus D_{n-1}),$$
where 
$$F_n=\left[
\begin{array}{cccccc}
1        &  b_0 & c_{0,1} & c_{0,2} & \ldots     & c_{0,n-1} \\
0        &  a_1 & b_1     & c_{1,2} &            & c_{1,n-1} \\  
0        &   0     & a_2     & b_2     &            &  \vdots      \\
0        &   0     &     0      & a_3     & \ddots     &              \\
\vdots   &         &            &   \ddots   & \ddots     &  b_{n-1}  \\
0        &         &  \ldots    &            & 0          &  a_n
\end{array}
\right].
$$
\end{lemma}
\begin{proof}
We have $D_0=1$ and then 
$$D_1=\left[
\begin{array}{cc}
1 & b_0 \\
0 & a_1 
\end{array}
\right]=F_1
\left[
\begin{array}{cc}
1 & 0 \\
0 & D_0 
\end{array}
\right].
$$
Assume the statement is true up to $n$. Then
$$D_{n+1}=\left[
\begin{array}{cc}
D_n & \\
    & J_{n+2,n+1}\ldots J_{2,1} \\
0   & 
\end{array}
\right].
$$ 
By the induction hypothesis, 
$$\left[
\begin{array}{c}
D_n \\
0 
\end{array}
\right]
=
\left[
\begin{array}{c}
F_n\left(1\oplus D_{n-1}\right) \\
0 
\end{array}
\right]
=
F_{n+1}\left[
\begin{array}{cc}
1   &   0     \\
0   & D_{n-1} \\
0   &   0
\end{array}
\right]
$$ 
and we notice that
$$J_{n+2,n+1}\ldots J_{2,1}=F_{n+1}\left[
\begin{array}{c}
0 \\
 J_{n+1,n}\ldots J_{2,1}
\end{array}
\right],
$$
so that 
$$D_{n+1}=F_{n+1}\left[
\begin{array}{ccc}
1   &  0      &  0 \\
0   & D_{n-1} &    \\ 
    &         &   J_{n+1,n}\ldots J_{2,1} \\
0   &  0      & 
\end{array}
\right]=F_{n+1}\left(1\oplus D_n\right).
$$
\end{proof}
The matrices $F_n$ have a very simple recursive multiplicative structure.
Actually it is convenient to make the dependence of $F_n$ on the 
Jacobi parameters more
explicit and introduce
$$F_{m,k}=\left[
\begin{array}{ccccc}
1    &   b_k   &  c_{k,k+1}  & \ldots  & c_{k,m-1}    \\
0    & a_{k+1} &    b_{k+1}  &         & c_{k+1,m-1}  \\
     &    0    &    a_{k+2}  &         &              \\
     &         &             & \ddots  & b_{m-1}      \\
     &         &             &    0    &   a_m 
\end{array}
\right]
$$
for $m\geq 1$ and $0\leq k<m-1$. In particular, $F_{n,0}=F_n$. 
We show that the building blocks of $F_{m,k}$ (consequently, of $D_n$)
are the $2\times 2$ matrices 
$$B_k=\left[
\begin{array}{cc}
 1  &   b_k \\
 0  &   a_{k+1}
\end{array}
\right], \quad k\geq 0,
$$
and 
the $(l-k+2)\times (l-k+2)$ matrices
$$C_{k,l}=\left[
\begin{array}{ccccc}
1    & 0    & \ldots   &  0         & c_{k,l} \\
0    & 1    &          &  0         &   0     \\     
     &      &          &  \ddots    &   0     \\
0    &      &          &  0         &   1
\end{array}
\right], \quad 0\leq k<l.
$$

\begin{lemma}\label{ldoua}
For $m\geq 1$ and $0\leq k<m-1$, 
$$F_{m,k}=\left(1\oplus F_{m,k+1}\right)G_{m,k},
$$
where 
$$G_{m,k}=C_{k,m-1}\left(C_{k,m-2}\oplus 1\right)
\ldots \left(C_{k,k+1}\oplus 1_{m-k-2}\right)\left(B_k\oplus 1_{m-k-1}\right).
$$
\end{lemma}
\begin{proof}
The proof is a straightforward calculation and can be omitted.
\end{proof}

Once again it is convenient to visualize all those matrix
multiplications by means of a transmission line
picture similar to the one in Figure ~3. Thus, 
Figure ~5 illustrates how to calculate $D_3$ by using 
Lemma ~\ref{luna} and Lemma ~\ref{ldoua}.
In particular, if $1$ is the imput at $A$ then at $B$ we read the 
expression of
$K(0,3)$ in terms of the Jacobi parameters, 
$$K(0,3)=b_0^3+b_0a_1c_{0,1}+a_1c_{0,1}b_0+a_1b_1c_{0,1}+a_1a_2c_{0,2}.
$$

Figure ~5 suggests a connection with weighted Lukasiewicz paths.
A Lukasiewicz path of length $n$ is a path in the positive 
quadrant of the lattice $\ZZ ^2$ which starts at $(0,0)$, 
ends at $(n,0)$, and consists of rise unit steps, 
horizontal unit steps, and fall steps of arbitrary depth. Let 
$\cL _{n,0}$ denote the set of Lukasiewicz paths of length $n$.

\begin{figure}[h]
\setlength{\unitlength}{3947sp}%
\begingroup\makeatletter\ifx\SetFigFont\undefined%
\gdef\SetFigFont#1#2#3#4#5{%
  \reset@font\fontsize{#1}{#2pt}%
  \fontfamily{#3}\fontseries{#4}\fontshape{#5}%
  \selectfont}%
\fi\endgroup%
\begin{picture}(4224,1230)(289,-886)
\thinlines
{\put(301,-661){\line( 1, 0){4200}}
}%
{\put(301,-361){\line( 1, 0){4200}}
}%
{\put(301,-61){\line( 1, 0){4200}}
}%
{\put(301,239){\line( 1, 0){4200}}
}%
\thicklines
{\put(901,-661){\line( 1, 1){300}}
}%
{\put(901,-661){\line( 1, 0){300}}
}%
{\put(1201,-361){\line( 1, 1){300}}
}%
{\put(1201,-361){\line( 1, 0){300}}
}%
{\put(1501,-661){\line( 1, 2){300}}
}%
{\put(1801,-661){\line( 1, 1){300}}
}%
{\put(1801,-661){\line( 1, 0){300}}
}%
{\put(2101,-61){\line( 1, 1){300}}
}%
{\put(2101,-61){\line( 1, 0){300}}
}%
{\put(2401,-361){\line( 1, 2){300}}
}%
{\put(2701,-661){\line( 1, 3){300}}
}%
{\put(3001,-361){\line( 1, 1){300}}
}%
{\put(3001,-361){\line( 1, 0){300}}
}%
{\put(3301,-661){\line( 1, 2){300}}
}%
{\put(3601,-661){\line( 1, 1){300}}
}%
{\put(3601,-661){\line( 1, 0){300}}
}%
\put(826,-511){$b_0$}%

\put(976,-826){$a_1$}%

\put(1876,-826){$a_2$}%

\put(3676,-826){$a_3$}%

\put(1276,-511){$a_1$}%

\put(1051,-211){$b_0$}%

\put(1726,-511){$b_1$}%

\put(2026,104){$b_0$}%

\put(2176,-211){$a_1$}%

\put(2926,-211){$b_1$}%

\put(3076,-511){$a_2$}%

\put(1451,-286){$c_{01}$}%

\put(2351, 14){$c_{01}$}%

\put(2551,-286){$c_{02}$}%

\put(3251,-286){$c_{12}$}%

\put(3501,-511){$b_2$}%

\put(376,-886){$A$}%

\put(4300,300){$B$}%

\end{picture}%

\caption{\mbox{ Transmission line representation for $D_3$ }}
\end{figure}

We also consider $\cL _{n,k}$ the set of paths of length 
$n$ in the positive quadrant, starting at $(0,0)$ and
consisting of the same type of steps as above, but ending at $(n,k)$.
We introduce a weigth on the elements of $\cL _{n,k}$ as 
follows. Let ${\bf p}\in \cL _{n,k}$
consists of steps ${\bf p}_1$, $\ldots $, ${\bf p}_n$.
Then 
$$w({\bf p})=\prod _{k=1}^nw({\bf p}_k)$$
and 
$$w({\bf p}_k)=\left\{
\begin{array}{cl}
a_l & \mbox{if ${\bf p}_k$ is a rise step $(j,l)\rightarrow (j+1,l+1)$
for some $j\geq 0$;} \\
b_l & \mbox{if ${\bf p}_k$ is a horizontal step $(j,l)\rightarrow (j+1,l)$
for some $j\geq 0$;} \\
c_{k,l} & \mbox{if ${\bf p}_k$ is a fall step $(j,l)\rightarrow (j+1,k)$
for some $j\geq l$.} 
\end{array}
\right.
$$

\begin{figure}[h]
\setlength{\unitlength}{3000sp}%
\begingroup\makeatletter\ifx\SetFigFont\undefined%
\gdef\SetFigFont#1#2#3#4#5{%
  \reset@font\fontsize{#1}{#2pt}%
  \fontfamily{#3}\fontseries{#4}\fontshape{#5}%
  \selectfont}%
\fi\endgroup%
\begin{picture}(1824,4147)(889,-3786)
\thicklines
{\put(1501,-61){\line( 1, 1){600}}
}%
{\put(1501,-61){\line( 1, -1){600}}
}%
{\put(1501,-61){\line( 1, 0){600}}
}%
{\put(1501,-61){\line( 1,-2){600}}
}%
{\put(1501,-661){\line( 1, 0){600}}
}%
{\put(1501,-661){\line( 1, 1){600}}
}%
{\put(1501,-661){\line( 1,-1){600}}
}%
{\put(1501,-2161){\line( 1, 0){600}}
}%
{\put(1501,-2161){\line( 1,-1){600}}
}%
{\put(1501,-2761){\line( 1, 1){600}}
}%
{\put(1501,-2761){\line( 1, 0){600}}
}%
{\put(1501,-2161){\line( 1,-2){600}}
}%
{\put(1501,-3361){\line( 1, 0){600}}
}%
{\put(1501,-3361){\line( 1, 1){600}}
}%
{\put(1501,-2761){\line( 1,-1){600}}
}%
{\put(1501,-2161){\line( 1, 1){300}}
}%
\thinlines
{\put(901,-3361){\line( 1, 0){1800}}
}%
\put(1576,314){$a_n$}%

\put(1651,-3586){$b_0$}%

\put(1426,-3136){$a_1$}%

\put(1351,-2461){$a_2$}%

\put(2101,-3061){$c_{01}$}%

\put(2101,-2536){$c_{12}$}%

\put(1801, 14){$b_{n-1}$}%

\put(851,-3786){$n-1$}%

\put(2376,-3786){$n$}%

\end{picture}%

\caption{\mbox{ Passing from paths of length $n-1$ to paths of lentgh $n$ }}
\end{figure}

\begin{theorem}\label{labella}
The Cholesky factor $D_n=\left[d_{i,j}\right]_{0\leq i,j\leq n}$ is given by
the formula
\begin{equation}\label{flabella}
d_{i,j}=\sum _{{\bf p}\in \cL _{j,i}}w({\bf p}), 
\quad i\leq j, \quad (i,j)\ne (0,0).
\end{equation}
\end{theorem}
\begin{proof}
We can prove the statement by induction on $n$. 
For $n\leq 3$, 
\eqref{flabella} is seen from Figure ~5. The general induction step
is provided by Lemma ~\ref{luna}. Thus, 
$$
\begin{array}{ccl}
d_{0,n}&=&b_0d_{0,n-1}+c_{0,1}d_{1,n-1}+\ldots +c_{0,n-1}d_{n-1,n-1} \\
d_{1,n}&=&a_1d_{0,n-1}+b_1d_{1,n-1}+\ldots +c_{1,n-1}d_{n-1,n-1} \\
       &\vdots & \\
d_{n,n}&=&a_nd_{n-1,n-1},
\end{array}
$$
and these relations are precisely those obtained by passing from weighted
paths of length $n-1$ to weighted paths of length $n$, 
as showed in Figure ~6.
\end{proof}

\bigskip
\noindent
{\bf Remarks}
$(a)$ For a Hankel kernel $K$, Theorem ~\ref{labella} 
reduces to well-known results in the 
combinatorial theory for orthogonal 
polynomials (on the real line), 
\cite{Fl}, \cite{Vi}. Indeed, by 
\eqref{carhank}
there are no fall steps of depth other than one. In this case,
the summation in \eqref{flabella} is only over Motzkin paths,
which is the classical formula in \cite{Fl}, \cite{Vi}.
It might be interesting to note that the correpsonding formula 
for orthogonal polynomials on the unit circle, 
\eqref{tsuma}, involves summation over labelled configurations 
in $\ZZ ^2$ rather than over weighted paths.

\smallskip
\noindent
$(b)$ There are other significant differences between the 
two parametrizations discussed in this paper. For instance, 
the parameters $\{\gamma _{k,j}\}$ have the 
following inheritance property: the parameters of the 
kernel $K^{(1)}=\left[K(l,m)\right]_{1\leq l,m}$ are precisely
$\{\gamma _{k,j}\}_{1\leq k<j}$. The Jacobi parameters do not
have such a property. Another difference involves computations of 
determinants. Thus, we have already notice the formula (we 
assume $K(0,0)=1$):
$$\det \left[K(l,m)\right]_{l,m=0}^n=
\prod _{k=1}^nf_k\prod _{0\leq i<j\leq n}d^2_{i,j},$$
and from Lemma ~\ref{luna} we deduce
$$\det \left[K(l,m)\right]_{l,m=0}^n=
\prod _{k=1}^na^{2(n-k)}_k,$$
which does not involve all the Jacobi parameters (up tu $n$). So the 
first determinant formula is much tighter in its parameters.

\smallskip
\noindent
$(c)$ Theorem ~\ref{labella} gives, in particular, that
\begin{equation}\label{ff}
K(0,n)=\sum _{{\bf p}\in \cL _{n,0}}w({\bf p}), \quad n\geq 1.
\end{equation}
Since the Jacobi parameters do not have the 
inheritance property mentioned above, 
we cannot have formulae of this type for any $K(l,n)$. Instead 
we have the following construction.

For $n$ fixed we consider admissible steps in the psitive quadrant of the 
following types: between the vertical lines $x=0$ and $x=n$ only
of Lukasiewicz type are allowed and between
the vertical lines $x=n$ and $x=2n$ only reflections with respect
to the line $x=n$ of Lukasiewicz steps are allowed (and they are
weighted with the complex conjugate of the weight of the 
reflected Lukasiewicz step, see Figure ~7).  
\begin{figure}[h]

\setlength{\unitlength}{3000sp}%
\begingroup\makeatletter\ifx\SetFigFont\undefined%
\gdef\SetFigFont#1#2#3#4#5{%
  \reset@font\fontsize{#1}{#2pt}%
  \fontfamily{#3}\fontseries{#4}\fontshape{#5}%
  \selectfont}%
\fi\endgroup%
\begin{picture}(2676,2413)(826,-2730)
\thinlines
{\put(901,-2161){\line( 1, 0){2400}}
}%
{\put(901,-2161){\line( 1, 1){1200}}
}%
{\put(2101,-961){\line( 1,-1){1200}}
}%
{\put(1501,-1561){\line( 1,-1){600}}
}%
{\put(2101,-2161){\line( 1, 1){600}}
}%
{\put(1501,-1561){\line( 1, 0){1200}}
}%
{\put(1501,-2161){\line( 1, 1){600}}
}%
{\put(2101,-1561){\line( 1,-1){600}}
}%
\thicklines
{\put(2101,-361){\line( 0,-1){2325}}
}%
\put(1126,-2461){$b_0$}%

\put(1651,-2461){$b_0$}%

\put(2326,-2461){$\overline{b}_0$}%

\put(2851,-2461){$\overline{b}_0$}%

\put(826,-1861){$a_1$}%

\put(1351,-1261){$a_2$}%

\put(2626,-1261){$a_2$}%

\put(3226,-1861){$a_1$}%

\put(1801,-1486){$b_1$}%

\put(2251,-1486){$\overline{b}_1$}%

\put(1501,-1786){$a_1$}%

\put(1726,-2011){$c_{01}$}%

\put(2251,-2011){$\overline{c}_{01}$}%

\put(2401,-1786){$a_1$}%

\end{picture}%

\caption{\mbox{ Admissible steps for $n=2$ }}
\end{figure}

Denote by $\cK _{n,k}$
the set of paths in the positive quadrant of $\ZZ ^2$ 
made of admissible steps, starting at $(0,0)$ and 
ending at $(n+k,0)$, $0\leq k\leq n$.
In particular, $\cK _{n,0}=\cL _{n,0}$. The weights are defined 
correspondingly. With these elements, 
we deduce from Theorem ~\ref{labella} that
\begin{equation}
K(l,n)=\sum _{{\bf p}\in \cK _{n,l}}w({\bf p}), \quad 0\leq l\leq n.
\end{equation}

\smallskip
\noindent
$(d)$ We notice that the proof of Theorem ~\ref{suma}
gives a formula for the Cholesky factor $D_n$ in terms
of Dyck type paths analogous to formula ~\eqref{flabella}.
\quad \quad \quad \quad \quad \quad \quad \quad \quad
$\Box $

\end{document}